\documentclass[11pt]{amsart}

\input{sommers_macros.tex}
 \usepackage{palatino}

\newcommand{\g}{{\mathfrak g}}
\newcommand{\tor}{{\mathfrak t}}

\newcommand{\uni}{{\mathfrak u}}

\newcommand{\orbit}{{\mathcal O}}

\newcommand{\al}{{\alpha}}
\newcommand{\coal}{{\alpha^{\scriptscriptstyle\vee}}}
\newcommand{\ro}{{\Phi}}

\newcommand{\complex}{{\mathbf C}}

\newcommand{\zz}{\mathbb Z}

\theoremstyle{definition}

\begin{document}

\setlength{\topmargin}{-.3in}
\setlength{\textheight}{9in} 
\setlength{\oddsidemargin}{-.2in}
\setlength{\evensidemargin}{-.2in} \setlength{\marginparsep}{0in}
\setlength{\marginparwidth}{0in}

\title {Normality of nilpotent varieties in $E_6$}

\author{ Eric Sommers}
\thanks{The author was supported in part by NSF grants DMS-0201826 and DMS-9729992.
He thanks the referee for a careful reading of the paper leading to its improvement.}

\address{
Department of Mathematics \\
University of Massachusetts--Amherst\\
Amherst, MA 01003 
}

\address{
Institute for Advanced Study \\
Princeton, NJ 08540
}

\date{4/15/03}

\email{esommers@math.umass.edu}

\begin{abstract}

We determine which nilpotent orbits in $E_6$ have closures which are normal
varieties and which do not.  
At the same time we are able to verify a conjecture in \cite{sommers:conjects} concerning
functions on non-special nilpotent orbits for $E_6$.

\end{abstract}
\maketitle

\section{Introduction}

The question of which nilpotent orbits in a simple Lie algebra 
(defined over the complex numbers) 
have normal closure has been studied by Kostant, Hesselink, Kraft-Procesi, 
Broer, and others.  
Kostant showed that the regular orbit has normal closure (that is, the
nilpotent cone is a normal variety) \cite{kostant:poly}.  
Kraft-Procesi showed that all 
nilpotent orbits in $\mathfrak{sl}_n(\complex)$ have normal closure~\cite{kraft-processi:normal}.
Vinberg-Popov showed
that the minimal orbit has normal closure \cite{vinberg-popov} and Hesselink
showed that several small orbits have normal closure~\cite{hesselink:normal}. 
Kraft-Procesi studied all nilpotent orbits in the classical groups and gave a method
to determine whether a nilpotent orbit has normal closure or not (their method does not
handle the very even orbits in the even orthogonal Lie algebras)~\cite{kraft-processi:normal2}.  
Kraft resolved the picture in $G_2$~\cite{kraft:normalG2} (see also 
Levasseur-Smith~\cite{levasseur-smith:prim}), 
and Broer resolved it in $F_4$~\cite{broer:normalF4}.  Broer
also showed that certain large orbits have normal closure 
(including the subregular orbit)~\cite{broer:vanishing}.
Over an algebraically closed field of good positive characteristic, Broer's work was extended by 
Kumar-Lauritzen-Thomsen and Thomsen \cite{klt:frob}, \cite{thomsen:normal}.
The methods in this paper constitute an extension of Broer's arguments 
in \cite{broer:vanishing} to smaller orbits (we will always work
over the complex numbers).  

Our main result is the determination of which orbits in $E_6$ have normal closure. 

\begin{thm}
The orbits in $E_6$ with normal closure are (in the notation of Bala-Carter)
$E_6$, $E_6(a_1)$, $D_5$, $E_6(a_3)$,
$D_5(a_1)$, $A_5$, $A_4+A_1$, $D_4$, $D_4(a_1)$, $2A_2+A_1$, 
$A_2+2A_1$, $A_2$, $3A_1$, $2A_1$, $A_1$, $0$.
The other $5$ orbits do not have normal closure.
\end{thm}

We also use the same techniques to prove directly 
a conjecture about functions on nilpotent orbit covers 
stated in \cite{sommers:conjects}.
Recall (after Broer) that a small representation is an irreducible highest weight
representation where twice a root is not a weight.  
Consider the following pairs of nilpotent orbits: $(A_5, E_6(a_3))$, 
$(2A_2+A_1, D_4(a_1))$, $(A_3+A_1, \tilde{D_4}(a_1))$, $(3A_1, A_2)$.
The first orbit in each pair is not special (in the sense of Lusztig)
and the second is its associated special orbit
(or in the case of $\tilde{D_4}(a_1)$, a 
$3$-fold cover of the associated special orbit that 
carries an action of the adjoint group of type $E_6$).
We show that 

\begin{thm}\label{small_conj}
The multiplicity of a small representation in the graded functions on the first orbit 
coincides with its multiplicity in the graded functions on the second orbit (or orbit cover). 
\end{thm}

Our proof is direct and realizes the functions on the first orbit as a quotient
of the functions on the second with a kernel that has no small representations in it.
In \cite{sommers:conjects} an analogous conjecture is stated and proved
for Springer fibers, but the proof is by calculating both sides and showing that the
multiplicities agree.  There does not appear at the present time to be an analog to the proof of
Theorem \ref{small_conj} for the Springer side of the picture.

\section{Notation}

Let $G$ be a simple, connected algebraic group defined over the complex numbers and 
$B$ a Borel subgroup containing a maximal torus $T$.
Let the character group of $T$ be $X^*(T)$ and $\ro$ the roots of $G$ with respect to $T$.

For any rational representation $\tau: B \to \mbox{GL}(V)$, 
let $V$ also denote the associated vector bundle $G \times^B V$ over $G/B$
when there is no ambiguity.  
In particular, if $\lambda \in X^*(T)=X^*(B)$, 
we write $\complex_{\lambda}$, or just $\lambda$, for the associated
line bundle on $G/B$ of the one-dimensional representation of $B$ coming from
$\lambda$.  We write $H^*(G/B, V)$ 
for the cohomology of $G/B$ in the sheaf of sections of $G \times^B V$.
If $P$ is a parabolic subgroup 
and $V$ is a representation of $P$, then
we can also consider the cohomology groups $H^*(G/P, V)$.  
Let $\complex[Y]$ denote the regular functions on an algebraic variety $Y$.

Let $\uni$ denote the Lie algebra of the unipotent radical of $B$.
We fix the {\it negative}
roots of $\ro$ to correspond to the weight spaces of $\uni$.  
Denote by $\ro^+$, $\ro^-$ the
positive and negative roots, respectively.  This choice also fixes
a set of simple roots $\SimpRoots=\{\al_i\}$.  Now let $W = N_G(T)/ T$ be the Weyl
group of $G$ and let $s_{\al}$ denote the reflection in 
the root $\al \in \ro$.  Let $\coal$ be the coroot for the root $\al \in \ro$ and
let $\langle \ \cdot \  , \ \cdot \ \rangle$ be the pairing of weights and coweights.  Denote by $P_{\al}$ the minimal parabolic subgroup containing $B$
corresponding to the simple root $\al$.

When $G$ is of type $E_6$, an element $\lambda \in X^*(T)$ 
of the form $\sum a_i \al_i$ where $\al_i \in \SimpRoots$ 
will be represented as
$\lbrace \begin{smallmatrix}
 a_1 & a_2 &  a_3 & a_4 & a_5 \\
   &   &  a_6 &  &
\end{smallmatrix} \rbrace$.  This also fixes our numbering of the simple roots.

We adopt the following notation for certain
subspaces of the nilradical $\uni$.
Let $h$ be an element of $\tor$, the Lie algebra of $T$.  
We can represent $h$ by the weighted Dynkin diagram 
with vertices labeled by $-\alpha_i(h)$ for the simple roots 
$\alpha_i \in \SimpRoots$.
We will denote the subspace $V = \oplus_{i \geq 2} \g_{i}$  
where $\g_i$ is the $i$-eigenspace of $\text{ad}(h)$ on $\g$
by putting brackets around the weighted Dynkin diagram for $h$.
Then $V$ will be a $B$-stable subspace of $\uni$ whenever 
all the vertices are labeled with non-negative real values.
For example, $\uni$ itself is represented by the diagram
$[\begin{smallmatrix}
2 & 2 &  2 & 2 & 2 \\
  &   &  2 &  &
\end{smallmatrix}]$.


\section{Method of Proof}

Assume that
$V \subset \uni$ is a subspace stable under the action of a parabolic subgroup $P$
which contains $B$.  Then $G \cdot V \subset \g$ (the $G$-saturation of $V$) is the
closure of a nilpotent orbit $\orbit$.
As explained in \cite{thomsen:normal}, the normality of the full nilpotent cone
implies that if the induced map $\complex[G \times^B \uni] \to \complex[G \times^B V]$
is surjective, then the closure of $\orbit$ is a normal variety.
Conversely, if $\bar{\orbit}$ is normal and the moment map 
$\mu: G \times^P V \to \bar{\orbit}$ is birational, 
then 
$\complex[G \times^B \uni] \to \complex[G \times^B V]$
is surjective.
The two key observations behind these statements are 
first, that $\complex[G \times^P V]= \complex[G \times^B V]$
and in fact more generally, for any $P$-representation $V$, that 
$H^i(G/B, V) = H^i(G/P, V)$ (see for example \cite{thomsen:normal});
and second, that when $\mu$ is birational 
$\complex[G \times^P V] = \complex[\orbit] = 
\complex[\bar{\orbit}]^{\text{norm}}$
(see for example \cite{kraft-processi:normal2}) where the latter
notation denotes the normalization.

Next consider the situation where $V_1 \subset V_2 \subset \uni$
and $V_i$ is stable under a parabolic subgroup $P_i$ which contains $B$.  
Let $i: G \times^B V_1 \to G \times^B V_2$ denote the inclusion.
Suppose that that the $G$-saturation of $V_2$ is known to be normal. 
Then it follows from the previous paragraph that we can deduce that 
the $G$-saturation of $V_1$ is normal if 
the induced map $i^*: \complex [G \times^B V_2] \to \complex[G \times^B V_1]$ is surjective
and the moment map $\mu: G \times^{P_2} V_2 \to G \cdot V_2$ is birational.
This will be our method of proof.  
We will also use the following elementary observation:  
if $i^*$ above is an isomorphism, then $G\cdot V_1 = G \cdot V_2$.  
This is an easy 
consequence of the fact that the moment maps are surjective.

In order to show that $i^*$ is surjective (respectively, an isomorphism),
we will start with the exact sequence of $B$-modules 
$$0 \to V_1 \to V_2 \to V_3 \to 0$$ (this defines $V_3$) and take the 
Koszul resolution of the dual sequence, obtaining the exact sequence of $B$-modules
\begin{equation} \label{koszul}  
\dots \to S^{n-j} V_2^* \otimes \wedge^{j} V_3^* \to \dots
\to S^{n-1} V_2^* \otimes V_3^*  \to S^n V_2^* \to S^n V_1^* \to 0.
\end{equation} 
Here, 
$S^n(-)$ denotes the $n$-th symmetric power
and $\wedge^j(-)$ the $j$-th exterior power.
By breaking the long exact sequence into short exact ones and taking the long exact
sequence in cohomology on $G/B$, we can often succeed in showing
that that induced map in cohomology $H^0(S^n V_2^*) \to H^0(S^n V_1^*)$
is surjective (resp, an isomorphism) for all $n \geq 0$.
This is sufficient to show that 
$i^*$ is surjective (resp, an isomorphism)
as we have the natural
isomorphism $\complex[G \times^B V]= \oplus_{n \geq 0} H^0(S^n V^*)$. 

\section{Tools}

We have three tools for showing that the induced map in 
cohomology $H^0(S^n V_2^*) \to H^0(S^n V_1^*)$
is surjective (resp, an isomorphism) for all $n \geq 0$.
Our first tool is the following key result of Demazure (see \cite{demazure:bott-simple}).

\begin{prop} \label{demazure}
Let $V$ be a rational representation of $B$ and assume that $V$ extends to 
a representation of the parabolic subgroup $P_{\alpha}$ where $\alpha$ is a simple root.
Let $\lambda \in X^*(T)$ 
be such that $m = \langle \lambda, \coal \rangle \geq -1$.
Then there is a $G$-module isomorphism 
$$ H^i( G/B, V \otimes {\lambda} ) = 
H^{i+1} ( G/B, V \otimes {\lambda \! - \!(m+1)\al} )
\text{ \ for \ all \ } i \in \zz.$$
In particular, if $m=-1$, then all cohomology groups vanish.
\end{prop}

Our second tool is a small
extension of a result of Broer (which relies on the vanishing theorem 
of Grauert and Reimenschneider) \cite{broer:vanishing}.  
Let $V$ be a subspace of $\uni$
stable under the action of a parabolic subgroup $P$ containing $B$
such that the moment map $\mu: G \times ^{P} V \to \g$ is 
generically finite.  This condition occurs in two special cases: 

{\bf Case 1.} 
$V=\uni_P$ is the Lie algebra of the unipotent radical of $P$. 

{\bf Case 2.} $V = \oplus_{i \geq 2} \g_{i}$ 
and $P$ is the parabolic subgroup with 
Lie algebra $\oplus_{i \geq 0} \g_{i}$ where $\g_i$ is the $i$-eigenspace for the semisimple
element of an $\mathfrak{sl}_2$-triple normalized so that $P$ contains $B$.  
Then,

\begin{prop} \label{broer}
For any dominant weight $\lambda \in X^*(P)$, we have
$$ H^i( G/P, S^n V^* \otimes {\omega} \otimes {\lambda} ) 
= H^i(G/B, S^n V^* \otimes {\omega} \otimes {\lambda})= 0 $$
for all $n \geq 0$ and $i > 0$,
where in case 1 above, $\omega = 0$, 
and in case 2 above, $\omega = \wedge^{\text top} \g_1$. 
\end{prop}

The proof for Case 1 is given in \cite{broer:vanishing} and the same proof 
also works for Case 2.

\medskip

Our third tool (which relies on the first tool) is proved in \cite{sommers:functions}.
Consider $G = SL_{l+1}(\complex)$.  
We label the simple roots $\Pi = \{ \al_1, \al_2, \dots, \al_l \}$ of $G$ 
so that consecutive indices are connected vertices in the Dynkin diagram
of type $A_l$.
Let $\{ \omega_j \}$ be the corresponding fundamental weights.

Let $P_m$ be the maximal (proper) parabolic subgroup of $G$ containing $B$
corresponding to all the simple roots except $\al_{m}$.  Denote
by $\uni_m$ the Lie algebra of the unipotent radical of $P_m$. 
The action of $P=P_m$ on $\uni_m$ gives
a representation of $P$ (and also $B$).  
Set $m' = \mbox{min}\{ m, l+1-m \}$.

\begin{prop} \label{sommers_A}
Let $r$ be an integer satisfying $2m'-2-l \leq r \leq 0$. 
Then there is a $G$-module isomorphism 
$$ H^i( G/B, S^{n} \uni^*_m \otimes r \omega_m )
= H^{i} ( G/B, S^{n + rm'} \uni^*_{l+1-m} \otimes -r \omega_{l+1-m})
\text{ \ for \ all \ } i, n \geq 0.$$
\end{prop}

This proposition can also often be applied in the more general setting where $G$ contains 
a Levi factor $L'$ of semisimple type $A_l$.  

More precisely, let $P$ be a parabolic subgroup of $G$ containing $B$ and let $L$ 
be the Levi factor of $P$ containing $T$.  Assume that 
$L$ contains simple factors of type $A_{m-1}$ and $A_{l-m}$ and 
these factors belong to a Levi subalgebra $L'$ of $G$ of type $A_l$.
Finally, assume that $[L, L'] \subset L'$.  
For ease of notation,
assume that the simple roots of $L'$ (which are also
simple roots for $G$) are labeled $\{ \al_1, \dots, \al_l \}$
and that $\al_m$ is not a simple root of $L$ (and hence $\al_i$ is a simple root
of $L$ for $i \neq m$).
The condition that $[L, L'] \subset L'$ is equivalent to saying that 
if a simple root of $G$ is connected in the Dynkin diagram to a simple root
of $L'$, then it is not a simple root of $L$.
Set $m' = \mbox{min}\{ m, l+1-m \}$.
Let $w_0$ denote the longest element of the Weyl group of $L'$.
Let $P_d$ denote the parabolic subgroup of $G$ containing $B$
with Levi factor equal to $L$
except we interchange the simple factors $A_{m-1}$ and $A_{l-m}$ in $L'$
(that is, we apply an outer automorphism to $L'$).
Let $\uni_P^*, \uni_{P_d}^*$ denote the Lie algebra of the unipotent radical of $P, P_d$,
respectively.

\begin{prop} \label{sommers}
Let $\lambda \in X^*(T)$ and set $r = \langle \lambda, \alpha_m \postcheck \rangle$.
Suppose that $\lambda$ satisfies 
$\langle \lambda, \alpha_i \postcheck \rangle = 0$  for $1 \leq i \neq m \leq l$ 
and $2m'-2-l \leq r \leq 0$. 

Then there is a $G$-module isomorphism 
$$ H^i( G/B, S^{n} \uni_{P}^* \otimes \lambda )
= H^{i} ( G/B, S^{n + rm'} \uni_{P_d}^* \otimes w_0(\lambda))
\text{ \ for \ all \ } i, n \geq 0.$$
\end{prop}

The proof in \cite{sommers:functions} also works for this more general case. 
Both of the previous two propositions show that when $n+rm' < 0$ all cohomology
groups vanish since the symmetric powers on the 
right-side of the equations are zero by definition.

\medskip

We now proceed through the nilpotent orbits of $E_6$ and determine whether
a given orbit has normal closure or not.
Since all cohomology considered henceforth will be the cohomology 
of vector bundles on $G/B$,
we omit the space $G/B$ from our notation for the cohomology.

\section{ $E_6, E_6(a_1), D_5, E_6(a_3)$}

$E_6$ has normal closure by Kostant \cite{kostant:poly} and the others by Broer \cite{broer:vanishing}.

\section{ $D_5(a_1)$}

The orbit $\orbit$ is Richardson for the 
parabolic subgroup whose Levi factor has semisimple part of type $A_2 + A_1$.
Therefore the closure of $\orbit$ is the $G$-saturation of 
$[\begin{smallmatrix}
0 & 0 &  2 & 2 & 0 \\
  &   &  2 &  &
\end{smallmatrix}]$.

\medskip
{\bf Step 1.}
Consider the short exact sequence of $B$-modules
$$0 \to [\begin{smallmatrix}
  0 & 1 &  1 & 2 & 0 \\
    &   &  2 &  &
\end{smallmatrix}]
\to 
[\begin{smallmatrix}
0 & 0 &  2 & 2 & 0 \\
   &   &  2 &  &
 \end{smallmatrix}] 
\to
V
\to
0$$
and take the Koszul resolution (equation \ref{koszul}) of the linear dual of this sequence.
There are only three terms in the resolution and the initial term equals  
$S^{n-1}[\begin{smallmatrix}
0 & 0 &  2 & 2 & 0 \\
   &   &  2 &  &
 \end{smallmatrix}]^* \otimes \complex_{\alpha_3}$
since $V = \complex_{-\alpha_3}$.
Since $\langle \alpha_3, \alpha_2 \postcheck \rangle = -1$
and the subspace $[\begin{smallmatrix}
 0 & 0 &  2 & 2 & 0 \\
    &   &  2 &  &
  \end{smallmatrix}]$ is stable under the parabolic $P_{\alpha_2}$,
Proposition \ref{demazure} with $m=-1$ implies that 

$$H^0(S^{n}[\begin{smallmatrix}
0 & 1 &  1 & 2 & 0 \\
   &   &  2 &  &
\end{smallmatrix}]^*) = 
H^0(S^{n}[\begin{smallmatrix}
0 & 0 &  2 & 2 & 0 \\
  &   &  2 &  &
\end{smallmatrix}]^*)
$$
for all $n$.
Consequently, $\bar{\orbit}$ equals $G \cdot [\begin{smallmatrix}
0 & 1 &  1 & 2 & 0 \\
   &   &  2 &  &
\end{smallmatrix}]$.

\medskip
{\bf Step 2.}
The closure of the orbit $E_6(a_3)$ is normal.
It is the $G$-saturation of 
$\uni_P = [\begin{smallmatrix}
 0 & 2 &  0 & 2 & 0 \\
   &   &  2 &  &
\end{smallmatrix}]$ with a birational moment map for the 
parabolic subgroup $P$ for which $\uni_P$ is the Lie algebra of its
unipotent radical.  The birationality follows (see \cite{mcgovern:regfn}) since this diagram
is the weighted Dynkin diagram for $E_6(a_3)$ and so 
for $e \in \uni$, the centralizer of $e$ in $G$ belongs
to $P$ (\cite{carter:book}, Proposition 5.7.1).

Consider the short exact sequence 
$$
0
\to
[\begin{smallmatrix}
0 & 1 &  1 & 2 & 0 \\
  &   &  2 &  &
\end{smallmatrix}]   
\to 
[\begin{smallmatrix}
0 & 2 &  0 & 2 & 0 \\
  &   &  2 &  &
\end{smallmatrix}]
\to 
V
\to
0.
$$

Taking the Koszul resolution of the dual sequence
yields
\begin{equation} \label{kossy}
0 \to S^{n-2} [ \begin{smallmatrix}
0 & 2 &  0 & 2 & 0 \\
  &   &  2 &  &
\end{smallmatrix} ]^*  \otimes 
\wedge^2 V^*
\to 
 S^{n-1}[\begin{smallmatrix}
0 & 2 &  0 & 2 & 0 \\
  &   &  2 &  &
\end{smallmatrix}]^*  \otimes 
V^*
\to 
S^{n}[\begin{smallmatrix}
0 & 2 &  0 & 2 & 0 \\
  &   &  2 &  &
\end{smallmatrix}]^*
\to 
S^{n}[\begin{smallmatrix}
0 & 1 &  1 & 2 & 0 \\
  &   &  2 &  &
\end{smallmatrix}]^* 
\to 0.
\end{equation}

\medskip
{\bf Step 3.}
We can simplify the cohomology of the two initial 
terms in the above exact sequence.
First, we compute that $\wedge^2 V^*$ equals $\complex_{\lambda}$ 
where $\lambda = \lbrace \begin{smallmatrix}
1 & 2 &  0 & 0 & 0 \\
  &   &  0 &  &
\end{smallmatrix} \rbrace$.
Applying Proposition \ref{demazure} 
with $m=0$ for the parabolic $P_{\alpha_3}$,
we get 
$$H^{i+1}(S^{n-2} [ \begin{smallmatrix}
0 & 2 &  0 & 2 & 0 \\
  &   &  2 &  &
\end{smallmatrix} ]^*  \otimes 
\wedge^2 V^*)
=
H^{i}(S^{n-2} [ \begin{smallmatrix}
 0 & 2 &  0 & 2 & 0 \\
   &   &  2 &  &
 \end{smallmatrix} ]^*  \otimes \lbrace \begin{smallmatrix}
 1 & 2 &  1 & 0 & 0 \\
   &   &  0 &  &
\end{smallmatrix} \rbrace)$$
for all $n \geq 2$ and all $i$.

Second, we can apply Proposition \ref{sommers} to the Levi factor
$L'$ with simple roots $\{ \al_3, \al_4, \al_5 \}$.
In this situation, $l=3$ and $m' = m = 2$.
The weight $\lbrace \begin{smallmatrix}
 1 & 2 &  1 & 0 & 0 \\
   &   &  0 &  &
\end{smallmatrix} \rbrace$
satisfies the hypothesis that 
its pairing is zero with $\al_3 \postcheck$ and $\al_5 \postcheck$
and that its pairing with $\al_4 \postcheck$ is $r=-1$.
In this case $P=P_d$ and we have
$$H^i(S^{n-2} [ \begin{smallmatrix}
0 & 2 &  0 & 2 & 0 \\
  &   &  2 &  &
\end{smallmatrix} ]^*  \otimes 
\lbrace \begin{smallmatrix}
 1 & 2 &  1 & 0 & 0 \\
   &   &  0 &  &
\end{smallmatrix} \rbrace) 
=
H^i(S^{n-4} [ \begin{smallmatrix}
0 & 2 &  0 & 2 & 0 \\
  &   &  2 &  &
\end{smallmatrix} ]^*  \otimes 
\lbrace \begin{smallmatrix}
 1 & 2 &  2 & 2 & 1 \\
   &   &  0 &  &
\end{smallmatrix} \rbrace).$$
We can then apply Proposition \ref{sommers} to the Levi factor 
$L'$ of type $A_2$ with simple roots $\{ \al_3,  \al_6 \}$.
Here $l=3$, $m=2$ and $m'=1$.
In our case, the weight 
$\lbrace \begin{smallmatrix}
 1 & 2 &  2 & 2 & 1 \\
   &   &  0 &  &
\end{smallmatrix} \rbrace$ satisfies the hypotheses with $r=-2$ since
its pairing with $\al_6 \postcheck$ is $-2$.
Therefore, 
$$H^i(S^{n-4} [ \begin{smallmatrix}
0 & 2 &  0 & 2 & 0 \\
  &   &  2 &  &
\end{smallmatrix} ]^*  \otimes 
\lbrace \begin{smallmatrix}
 1 & 2 &  2 & 2 & 1 \\
   &   &  0 &  &
\end{smallmatrix} \rbrace)
=H^i(S^{n-6} [ \begin{smallmatrix}
0 & 2 &  2 & 2 & 0 \\
  &   &  0 &  &
\end{smallmatrix} ]^*  \otimes 
\lbrace \begin{smallmatrix}
 1 & 2 &  4 & 2 & 1 \\
   &   &  2 &  &
\end{smallmatrix} \rbrace).$$
Two more applications of the proposition (to the symmetric $A_2$ factors on either end;
in both situations $m'=1$ and $r=-1$)
yield that the latter is isomorphic to 
$$H^i(S^{n-8} [ \begin{smallmatrix}
 2 & 0 &  2 & 0 & 2 \\
   &   &  0 &  &
 \end{smallmatrix} ]^*  \otimes \lbrace \begin{smallmatrix}
 2 & 3 &  4 & 3 & 2 \\
   &   &  2 &  &
\end{smallmatrix} \rbrace).$$

Then by Proposition \ref{broer} these cohomologies vanish for $i>0$.
 
Finally, we can show that 
$H^i(S^{n-1}[\begin{smallmatrix}
0 & 2 &  0 & 2 & 0 \\
  &   &  2 &  &
\end{smallmatrix}]^*  \otimes 
V^*) = 0$
for all $i \geq 0$.
This is because we have the short exact sequence of $B$-modules 
$0 \to \complex_{\lbrace \begin{smallmatrix}
0 & 1 &  0 & 0 & 0 \\
  &   &  0 &  &
\end{smallmatrix} \rbrace}
\to  V^* \to
\complex_{\lbrace \begin{smallmatrix}
1 & 1 &  0 & 0 & 0 \\
  &   &  0 &  &
\end{smallmatrix} \rbrace}
 \to 0$.
Since both the weights on either end of the sequence have inner product $m=-1$ 
with $\alpha_3 \postcheck$ and 
$[\begin{smallmatrix}
0 & 2 &  0 & 2 & 0 \\
  &   &  2 &  &
\end{smallmatrix}]$ is $P_{\alpha_3}$-stable, 
we get the vanishing by Proposition \ref{demazure}. 

\medskip
{\bf Step 4.}
Breaking the Koszul sequence (\ref{kossy}) into (two) 
short exact sequences, taking the long exact sequence in cohomology with 
respect to $G/B$, and using the results in Step 3 above,
yields the following exact sequence for all $n$

$$0 \to
H^0 (S^{n-8} [ \begin{smallmatrix}
2 & 0 &  2 & 0 & 2 \\
  &   &  0 &  &
\end{smallmatrix} ]^*  \otimes \lbrace \begin{smallmatrix}
2 & 3 &  4 & 3 & 2 \\
  &   &  2 &  &
\end{smallmatrix} \rbrace)
\to 
H^0(S^{n}[\begin{smallmatrix}
0 & 2 &  0 & 2 & 0 \\
  &   &  2 &  &
\end{smallmatrix}]^*) 
\to 
H^0(S^{n}[\begin{smallmatrix}
0 & 1 &  1 & 2 & 0 \\
  &   &  2 &  &
\end{smallmatrix}]^*) 
\to 0,
$$
and thus the closure of $\orbit$ is normal.

\section{ $A_5$} \label{section_a5}

The closure of the orbit $\orbit$ of type $A_5$ 
is the $G$-saturation of 
$[\begin{smallmatrix}
2 & 1 &  0 & 1 & 2 \\
  &   &  1 &  &
\end{smallmatrix}]$ (this is its weighted Dynkin diagram).  
As in the previous section we utilize the normality of the closure 
of the orbit $E_6(a_3)$.

Consider the short exact sequence
$$
0
\to
[\begin{smallmatrix}
2 & 1 &  0 & 1 & 2 \\
  &   &  1 &  &
\end{smallmatrix}]
\to
[\begin{smallmatrix}
2 & 0 &  2 & 0 & 2 \\
  &   &  0 &  &
\end{smallmatrix}]
\to
V
\to
0
$$
and take the Koszul resolution of its dual
(since the dimension of $V$ is four, the resolution has six terms). 
We can simplify the cohomology of the four initial terms in the resolution.
First, 
 we can show that
$$H^i(S^{n-j} [ \begin{smallmatrix}
2 &  0 &  2 & 0 & 2 \\
  &   &   0 &  &
\end{smallmatrix} ]^*  \otimes 
\wedge^j V^*)=0$$
for all $i$,$n$ and for $j = 1, 3$.  This is because 
$\wedge^j V^*$ for $j=1,3$ can be filtered by 
$B$-subrepresentations whose quotients
are one-dimensional and which have total vanishing cohomology by 
Proposition \ref{demazure} as in the last part of Step 3 of the previous orbit.
However, the case of $j=2$ is more difficult.  To deduce that we have total vanishing
for $j=2$ requires the study of a specific six-dimensional bundle on the product of three
projective lines.

\begin{lem} \label{a1_thrice}
Let $G = SL_2 \times SL_2 \times SL_2$ and let $B$ be
a Borel subgroup of $G$.  Let $U$ be the eight-dimensional irreducible
representation of $G$ which is the tensor product of 
the three standard representations for each $SL_2$ factor of $G$.
Let $U'$ be the four-dimensional $B$-stable subspace of $U$
containing the four lowest weight spaces of $U$.  
Then $H^i(G/B, \wedge^2 U') = 0 $ for all $i \geq 0$.
\end{lem}

\begin{proof}
We study the Koszul resolution
$$0 \to  \wedge^2 U' \to U \otimes U' \to S^2 U \to S^2(U/U') \to 0.$$

Let $\al_1, \al_2, \al_3$ be the three simple roots of $G$.
The restriction of $U'$ to a maximal torus $T$ in $B$ 
yields the weights  $(-1,-1,-1), (1,-1,-1), 
(-1,1,-1),(-1, -1,1)$ with respect to the simple roots.
Proposition \ref{demazure} thus implies that 
$H^i(G/B, U \otimes U') = 0 $ for all $i\geq 0$
(this is identical to the proof used above for $V^*$ and $j=1$).
We deduce that the kernel of the map from $H^0(G/B, S^2 U)$ to $H^0(G/B,  S^2(U/U'))$
is isomorphic to $H^1(G/B, \wedge^2 U')$.

Now $\wedge^2 U'$ can be filtered by $B$-subrepresentations 
with bases
$$ \{ (-2,-2,0) \} \subset \{ (-2,-2,0), (-2, 0,-2), (-2, 0, 0) \} \subset $$
$$\{ (-2,-2,0), (-2, 0,-2), (-2, 0, 0), (0, -2, -2), (0, -2, 0) \} \subset 
\wedge^2 U' ,$$
where we use the weight to denote the corresponding weight vector.
The quotients of these subrepresentations yield two line bundles
and two two-dimensional vector bundles.  The latter have total vanishing
by Proposition \ref{demazure} (applied to the case where the $V$ in that proposition
is a two-dimensional irreducible representation and $m=-1$).
The line bundles, however, do not have total vanishing.  Indeed 
the line bundle of weight $(-2,-2,0)$
has $H^2$ equal to the trivial representation of $G$;
and the line bundle of weight $(0,0,-2)$ has
$H^1$ equal to the trivial representation of $G$
(this follows from Bott-Borel-Weil, which in this setting is equivalent
to Proposition \ref{demazure}).
Consequently, we find that either  
$H^i(G/B, \wedge^2 U')$ vanishes for all $i$ or
it is equal to the trivial representation for $i=1,2$ and vanishes otherwise.
In the latter case, it would follow that $H^0(G/B, S^2 U)$ contains a copy of
the trivial representation.

But a calculation shows that 
$S^2 U$ is the direct sum of three irreducible $3$-dimensional representations of $G$ and one 
irreducible $27$-dimensional representation.  Since it does not contain a copy
of the trivial representation and since $H^0(G/B, S^2 U) = S^2 U$ as $U$ is a
representation of $G$, the latter scenario can not occur. 
\end{proof}

We now apply this lemma to our situation.
The $B$-representation 
$$S^{n-j} [ \begin{smallmatrix}
2 &  0 &  2 & 0 & 2 \\
   &   &   0 &  &
\end{smallmatrix} ]^*$$ extends to a $P$-representation where $P$ is a parabolic
containing $B$ such that $P/B$ is isomorphic to the product of three projective lines. 
The representation $V^*$ then yields a bundle on $P/B$ isomorphic to 
the bundle determined by $U'$ from the lemma.  Hence the lemma and 
a spectral sequence argument as in the proof of Proposition \ref{demazure}
in \cite{demazure:bott-simple} yields the result.

\medskip

Next, 
we compute that $\wedge^4(V^*) = \complex_{\lambda}$ where 
$\lambda = \lbrace \begin{smallmatrix}
0 & 1 &  4 & 1 & 0 \\
  &   &  1 &  &
\end{smallmatrix} \rbrace$ in the basis of simple roots.
Then after three applications of Proposition \ref{demazure}
with $m=0$ for each of the parabolics $P_{\alpha_2}$, 
$P_{\alpha_4}$, and $P_{\alpha_6}$, we get
$$H^{i+3}(S^{n-4} [ \begin{smallmatrix}
2 &  0 &  2 & 0 & 2 \\
  &   &   0 &  &
\end{smallmatrix} ]^*  \otimes 
\wedge^4 V^*)
=
H^{i}(S^{n-4} [ \begin{smallmatrix}
2 &  0 &  2 & 0 & 2 \\
  &   &   0 &  &
\end{smallmatrix} ]^*  \otimes 
\lbrace \begin{smallmatrix}
0 & 2 &  4 & 2 & 0 \\
  &   &  2 &  &
\end{smallmatrix} \rbrace).
$$
Now we use Proposition \ref{sommers}  
three times (to each of the extreme $A_2$ factors)
and get 
$$H^{i}(S^{n-4} [ \begin{smallmatrix}
2 &  0 &  2 & 0 & 2 \\
  &   &   0 &  &
\end{smallmatrix} ]^*  \otimes 
\lbrace \begin{smallmatrix}
0 & 2 &  4 & 2 & 0 \\
  &   &  2 &  &
\end{smallmatrix} \rbrace)
=
H^{i}(S^{n-8} [ \begin{smallmatrix}
0 &  2 &  2 & 2 & 0 \\
  &   &   0 &  &
\end{smallmatrix} ]^*  \otimes 
\lbrace \begin{smallmatrix}
2 & 4 &  4 & 4 & 2 \\
  &   &  2 &  &
\end{smallmatrix} \rbrace)
=$$
$$=H^{i}(S^{n-10} [ \begin{smallmatrix}
0 &  2 &  0 & 2 & 0 \\
  &   &   2 &  &
\end{smallmatrix} ]^*  \otimes 
\lbrace \begin{smallmatrix}
2 & 4 &  6 & 4 & 2 \\
  &   &  4 &  &
\end{smallmatrix} \rbrace)
.$$
Then Proposition \ref{broer} implies that the cohomology of the latter vector bundle
is trivial if $i>0$ for all $n$.

Now we can finish by breaking the Koszul resolution into short exact sequences
and taking the long exact sequence in cohomology.
We thus have that
$$0 \to H^0(S^{n-10} [ \begin{smallmatrix}
0 &  2 &  0 & 2 & 0 \\
  &   &   2 &  &
\end{smallmatrix} ]^*  \otimes \lbrace \begin{smallmatrix}
2 & 4 &  6 & 4 & 2 \\
  &   &  4 &  &
\end{smallmatrix} \rbrace) 
\to 
H^0(S^{n}[\begin{smallmatrix}
2 & 0 &  2 & 0 & 2 \\
  &   &  0 &  &
\end{smallmatrix}]^*) 
\to 
H^0(S^{n}[\begin{smallmatrix}
2 & 1 &  0 & 1 & 2 \\
  &   &  1 &  &
\end{smallmatrix}]^*) 
\to 0$$
is exact, proving the normality of the closure of $\orbit$.

\section{ $A_4 + A_1$}

This orbit is Richardson for any parabolic subgroup whose Levi
subgroup has semisimple part of type $A_2+ 2A_1$. 
Hence its closure equals $G \cdot [\begin{smallmatrix}
0 & 0 &  2 & 2 & 0 \\
  &   &  0 &  &
\end{smallmatrix}]$.
We prove normality by using the (just proved) normality of $D_5(a_1)$.
The closure of $D_5(a_1)$ equals
the $G$-saturation of 
$[\begin{smallmatrix}
 0 & 0 &  2 & 0 & 2 \\
   &   &  2 &  &
 \end{smallmatrix}]$ 
with birational moment map for the maximal parabolic $P$ which stabilizes this subspace. 
The birationality follows since the centralizer in $G$ of any element 
in $D_5(a_1)$ is connected and hence the centralizer in $G$ of $e \in
[\begin{smallmatrix}
 0 & 0 &  2 & 0 & 2 \\
   &   &  2 &  &
 \end{smallmatrix}]$ equals the centralizer of $e$ in $P$ (see \cite{carter:book}).

Consider the short exact sequence
$$
0
\to
[\begin{smallmatrix}
0 & 0 &  2 & 2 & 0 \\
  &   &  0 &  &
\end{smallmatrix}]
\to
[\begin{smallmatrix}
0 & 0 &  2 & 2 & 0 \\
  &   &  2 &  &
\end{smallmatrix}]
\to
V
\to
0
$$
and take the Koszul resolution of its dual (there are only three terms).

We have $V^* = \complex_{\lambda}$ where 
$\lambda = \lbrace \begin{smallmatrix}
0 & 0 &  0 & 0 & 0 \\
  &   &  1 &  &
\end{smallmatrix} \rbrace$.
Now we use Proposition \ref{sommers} three times to get 
$$H^{i}(S^{n-1} [ \begin{smallmatrix}
0 &  0 &  2 & 2 & 0 \\
  &   &   2 &  &
\end{smallmatrix} ]^*  \otimes 
\lbrace \begin{smallmatrix}
0 & 0 &  0 & 0 & 0 \\
  &   &  1 &  &
\end{smallmatrix} \rbrace)
=
H^{i}(S^{n-2} [ \begin{smallmatrix}
2 &  0 &  0 & 2 & 0 \\
  &   &   2 &  &
\end{smallmatrix} ]^*  \otimes 
\lbrace \begin{smallmatrix}
1 & 1 &  1 & 0 & 0 \\
  &   &  1 &  &
\end{smallmatrix} \rbrace)
= $$
$$H^{i}(S^{n-4} [ \begin{smallmatrix}
2 &  0 &  2 & 0 & 0 \\
  &   &   2 &  &
\end{smallmatrix} ]^*  \otimes 
\lbrace \begin{smallmatrix}
1 & 2 &  3 & 2 & 1 \\
  &   &  1 &  &
\end{smallmatrix} \rbrace)
=
H^{i}(S^{n-5} [ \begin{smallmatrix}
2 &  0 &  2 & 0 & 0 \\
  &   &   2 &  &
\end{smallmatrix} ]^*  \otimes 
\lbrace \begin{smallmatrix}
1 & 2  &  3 & 2 & 1\\
  &   &  2 &  &
\end{smallmatrix} \rbrace)
.$$
Then Proposition \ref{broer} implies that the cohomology of the latter vector bundle
is trivial if $i>0$ for all $n$.

We thus have a short exact sequence in cohomology
$$0 \to H^0(S^{n-5} [ \begin{smallmatrix}
2 &  0 &  2 & 0 & 0 \\
  &   &   2 &  &
\end{smallmatrix} ]^*  \otimes 
\lbrace \begin{smallmatrix}
1 & 2  &  3 & 2 & 1\\
  &   &  2 &  &
\end{smallmatrix} \rbrace)
\to 
H^0(S^{n}[\begin{smallmatrix}
0 &  0 &  2 & 2 & 0 \\
  &   &   2 &  &
\end{smallmatrix}]^*) 
\to 
H^0(S^{n}[\begin{smallmatrix}
0 & 0 &  2 & 2 & 0 \\
  &   &  0 &  &
\end{smallmatrix}]^*) 
\to 0,$$
proving normality.

\section{ $D_4$}

The closure of $D_4$ equals $G \cdot [\begin{smallmatrix}
0 & 0 &  2 & 0 & 0 \\
  &   &  2 &  &
\end{smallmatrix}]$.   We prove normality by 
again using the fact that the closure of $D_5(a_1)$ is normal.

Hence we study
$$
0
\to
[\begin{smallmatrix}
0 & 0 &  2 & 0 & 0 \\
  &   &  2 &  &
\end{smallmatrix}]
\to
[\begin{smallmatrix}
0 & 0 &  2 & 0 & 2 \\
  &   &  2 &  &
\end{smallmatrix}]
\to
V
\to
0
$$
and take the Koszul resolution of its dual (there are four terms).

For the first term of the resolution, 
$$H^i(S^{n-2} [ \begin{smallmatrix}
0 &  0 &  2 & 0 & 2 \\
  &   &   2 &  &
\end{smallmatrix} ]^*  \otimes \wedge^2 V^*)
 = H^i(S^{n-2} [ \begin{smallmatrix}
0 &  0 &  2 & 0 & 2 \\
  &   &   2 &  &
\end{smallmatrix} ]^*  \otimes \lbrace \begin{smallmatrix}
0 & 0 &  0 & 1 & 2 \\
  &   &  0 &  &
\end{smallmatrix} \rbrace),$$
and then by Proposition \ref{sommers},
$$ H^i(S^{n-2} [ \begin{smallmatrix}
0 &  0 &  2 & 0 & 2 \\
  &   &   2 &  &
\end{smallmatrix} ]^*  \otimes \lbrace \begin{smallmatrix}
0 & 0 &  0 & 1 & 2 \\
  &   &  0 &  &
\end{smallmatrix} \rbrace) 
=
H^i(S^{n-4} [ \begin{smallmatrix}
0 &  2 &  0 & 0 & 2 \\
  &   &   2 &  &
\end{smallmatrix} ]^*  \otimes \lbrace \begin{smallmatrix}
1 & 2 &  2 & 2 & 2 \\
  &   &  0 &  &
\end{smallmatrix} \rbrace) 
=$$
$$H^i(S^{n-6} [ \begin{smallmatrix}
0 &  2 &  0 & 2 & 2 \\
  &   &   0 &  &
\end{smallmatrix} ]^*  \otimes \lbrace \begin{smallmatrix}
1 & 2 &  4 & 4 & 2 \\
  &   &  2 &  &
\end{smallmatrix} \rbrace)
=
H^i(S^{n-8} [ \begin{smallmatrix}
0 & 0 &  2 & 2 & 2 \\
  &   &  0 &  &
\end{smallmatrix} ]^*  \otimes \lbrace \begin{smallmatrix}
2 & 4 &  6 & 4 & 2 \\
  &   &  3 &  &
\end{smallmatrix} \rbrace),
$$
and thus all these groups vanish for $i>0$
by Proposition \ref{broer}.

On the other hand, Proposition \ref{demazure},
with $m=-1$ and $P=P_{\alpha_4}$, gives that
$H^i(S^{n-1} [ \begin{smallmatrix}
0 &  0 &  2 & 0 & 2 \\
  &   &   2 &  &
\end{smallmatrix} ]^*  \otimes \lbrace \begin{smallmatrix}
0 & 0 &  0 & 0 & 1 \\
  &   &  0 &  &
\end{smallmatrix} \rbrace) = 0$ for all $i$, and so 
$$H^i(S^{n-1} [ \begin{smallmatrix}
0 &  0 &  2 & 0 & 2 \\
  &   &   2 &  &
\end{smallmatrix} ]^* \otimes V^*) = 
H^i(S^{n-1} [ \begin{smallmatrix}
0 &  0 &  2 & 0 & 2 \\
  &   &   2 &  &
\end{smallmatrix} ]^*  \otimes \lbrace \begin{smallmatrix}
0 & 0 &  0 & 1 & 1 \\
  &   &  0 &  &
\end{smallmatrix} \rbrace).$$

Now consider
the exact sequence 
\begin{equation} \label{kos2:a_4+a_1}
\dots \to H^i(S^{n-2} [ \begin{smallmatrix}
0 &  0 &  2 & 2 & 2 \\
  &   &   2  &  &
\end{smallmatrix} ]^*  \otimes \lbrace \begin{smallmatrix}
0 & 0 &  0 & 2 & 1 \\
  &   &  0 &  &
\end{smallmatrix} \rbrace) \to
\end{equation}
$$
H^i(S^{n-1}[\begin{smallmatrix}
0 & 0 &  2 & 2 & 2 \\
  &   &  2  &  &
\end{smallmatrix}]^* \otimes \lbrace \begin{smallmatrix}
0 & 0 &  0 & 1 & 1 \\
  &   &  0 &  &
\end{smallmatrix} \rbrace) 
\to 
H^i(S^{n-1} [ \begin{smallmatrix}
0 &  0 &  2 & 0 & 2 \\
  &   &   2 &  &
\end{smallmatrix} ]^*  \otimes \lbrace \begin{smallmatrix}
0 & 0 &  0 & 1 & 1 \\
  &   &  0 &  &
\end{smallmatrix} \rbrace)
\to \dots$$
obtained from the obvious three term Koszul resolution tensored with
the weight 
$\lbrace \begin{smallmatrix}
0 & 0 &  0 & 1 & 1 \\
  &   &  0 &  &
\end{smallmatrix} \rbrace$.
Now
$$H^i(S^{n-1}[\begin{smallmatrix}
0 & 0 &  2 & 2 & 2 \\
  &   &  2  &  &
\end{smallmatrix}]^* \otimes \lbrace \begin{smallmatrix}
0 & 0 &  0 & 1 & 1 \\
  &   &  0 &  &
\end{smallmatrix} \rbrace) =
H^i(S^{n-4}[\begin{smallmatrix}
2 & 2 &  0 & 0 & 2 \\
  &   &  2  &  &
\end{smallmatrix}]^* \otimes \lbrace \begin{smallmatrix}
1 & 2 &  3 & 2 & 1 \\
  &   &  2 &  &
\end{smallmatrix} \rbrace),$$ which vanishes for $i>0$,
and
$$H^i(S^{n-2}[\begin{smallmatrix}
0 & 0 &  2 & 2 & 2 \\
  &   &  2  &  &
\end{smallmatrix}]^* \otimes \lbrace \begin{smallmatrix}
0 & 0 &  0 & 2 & 1 \\
  &   &  0 &  &
\end{smallmatrix} \rbrace)=
H^i(S^{n-8}[\begin{smallmatrix}
2 & 2 &  2 & 0 & 0 \\
  &   &  2  &  &
\end{smallmatrix}]^* \otimes \lbrace \begin{smallmatrix}
2 & 4 &  6 & 4 & 2 \\
  &   &  3 &  &
\end{smallmatrix} \rbrace),$$ which also vanishes for $i>0$.
Therefore,
$H^i(S^{n-1} [ \begin{smallmatrix}
0 &  0 &  2 & 0 & 2 \\
  &   &   2 &  &
\end{smallmatrix} ]^*  \otimes \lbrace \begin{smallmatrix}
0 & 0 &  0 & 1 & 1 \\
  &   &  0 &  &
\end{smallmatrix} \rbrace)$ vanishes for $i>0$ by Equation \ref{kos2:a_4+a_1}.
Finally,
going back to the original Koszul resolution, breaking it into two 
short exact sequences, and taking cohomology
yields 
that $H^0(S^{n}[\begin{smallmatrix}
0 & 0 &  2 & 0 & 2 \\
  &   &  2 &  &
\end{smallmatrix}]^*)$ 
surjects onto 
$H^0(S^{n}[\begin{smallmatrix}
0 & 0 &  2 & 0 & 0 \\
  &   &  2 &  &
\end{smallmatrix}]^*)$, 
proving normality. 



\section{ $D_4(a_1)$}

There is an exact sequence
\begin{equation} \label{kos1:d_4(a_1)}
\dots \to H^i(S^{n-1} [ \begin{smallmatrix}
0 &  0 &  2 & 0 & 0 \\
  &   &   2 &  &
\end{smallmatrix} ]^*  \otimes \lbrace \begin{smallmatrix}
0 & 0 &  0 & 0 & 0 \\
  &   &  1 &  &
\end{smallmatrix} \rbrace) 
\to 
H^i(S^{n}[\begin{smallmatrix}
0 & 0 &  2 & 0 & 0 \\
  &   &  2 &  &
\end{smallmatrix}]^*) 
\to 
H^i(S^{n}[\begin{smallmatrix}
0 & 0 &  2 & 0 & 0 \\
  &   &  0 &  &
\end{smallmatrix}]^*) 
\to \dots.
\end{equation}
Then
$H^i(S^{n-1} [ \begin{smallmatrix}
0 &  0 &  2 & 0 & 0 \\
  &   &   2 &  &
\end{smallmatrix} ]^*  \otimes \lbrace \begin{smallmatrix}
0 & 0 &  0 & 0 & 0 \\
  &   &  1 &  &
\end{smallmatrix} \rbrace)$
equals
$H^i(S^{n-5} [ \begin{smallmatrix}
0 &  0 &  2 & 0 & 0 \\
  &   &   2 &  &
\end{smallmatrix} ]^*  \otimes \lbrace \begin{smallmatrix}
1 & 2 &  3 & 2 & 1 \\
  &   &  2 &  &
\end{smallmatrix} \rbrace)$, which vanishes for $i>0$.
The desired surjection follows as does normality since $D_4$ is
normal and the appropriate moment map is birational.

\section{ $A_3 + A_1$}

Although the closure of $\orbit$ is not normal, 
the regular functions on $\orbit$ are 
a quotient of the functions on a $3$-fold cover of the orbit $D_4(a_1)$ and 
we will use this in the proof of Theorem \ref{small_conj}.  

The weighted Dynkin diagram of $A_3 +A_1$ is
$ \begin{smallmatrix}
0 &  1 &  0 & 1 & 0 \\
  &   &   1 &  &
\end{smallmatrix}$ and thus
$\complex[\orbit] = \complex[G \times^B [ \begin{smallmatrix}
0 &  1 &  0 & 1 & 0 \\
  &   &   1 &  &
\end{smallmatrix} ]]$.
On the other hand, the subspace $[ \begin{smallmatrix}
0 &  0 &  0 & 2 & 0 \\
  &   &   2 &  &
\end{smallmatrix} ]$ has $G$-saturation $D_4(a_1)$ and the moment map (for the
natural parabolic subgroup $P$ coming from the zeros in the diagram) is generically
3-to-1.  This is a consequence of the fact that the centralizer in $P$ of a Richardson 
element $e$ has index three in the centralizer in $G$ of $e$.
It follows that   
$\complex[G \times^B [ \begin{smallmatrix}
0 &  0 &  0 & 2 & 0 \\
  &   &   2 &  &
\end{smallmatrix} ]]$
equals the functions on a $3$-fold cover of $D_4(a_1)$.

Consider the subspace $$U = [ \begin{smallmatrix}
0 &  0 &  2 & 0 & 0 \\
  &   &   0 &  &
\end{smallmatrix} ] \cap 
[ \begin{smallmatrix}
0 &  0 &  0 & 2 & 0 \\
  &   &   2 &  &
\end{smallmatrix} ]$$ and the two exact sequences
\begin{equation} \label{a4,1}
0 \to
U \to
[ \begin{smallmatrix}
0 &  0 &  2 & 0 & 0 \\
  &   &   0 &  &
\end{smallmatrix} ] \to
V_1 \to 0
\end{equation} 
and 
\begin{equation} \label{a4,2}
0 \to
U \to
[ \begin{smallmatrix}
0 &  0 &  0 & 2 & 0 \\
  &   &   2 &  &
\end{smallmatrix} ] \to
V_2 \to 0.
\end{equation}
Analysis of the Koszul resolution of the dual of Equation \ref{a4,2} 
and several applications of Proposition \ref{demazure}
shows that 
$H^i(S^n U^*) = H^i(S^n[ \begin{smallmatrix}
0 &  0 &  0 & 2 & 0 \\
  &   &   2 &  &
\end{smallmatrix} ]^*)$ and thus $\complex[G \times^B U]$  also equals 
the functions on a $3$-fold cover of $D_4(a_1)$.  
We omit the details as they are similar to those in the previous sections.

Next consider the exact sequence
$$0 \to
[ \begin{smallmatrix}
0 &  1 &  0 & 1 & 0 \\
  &   &   1 &  &
\end{smallmatrix} ]
 \to
U \to
V_3 \to 0.$$
This leads to the exact sequence
$$\dots \to H^i(S^{n-3} U^*
\otimes \lbrace \begin{smallmatrix}
1 & 2 &  3 & 2 & 1 \\
  &   &  1 &  &
\end{smallmatrix} \rbrace) 
\to 
H^i(S^{n} U^*)
\to 
H^i(S^{n}[ \begin{smallmatrix}
0 &  1 &  0 & 1 & 0 \\
  &   &   1 &  &
\end{smallmatrix} ]^*) 
\to \dots.$$
Then taking the Koszul resolution of the dual of Equation \ref{a4,1} and
tensoring with $\lbrace \begin{smallmatrix}
1 & 2 &  3 & 2 & 1 \\
  &   &  1 &  &
\end{smallmatrix} \rbrace$ yields, after some work, 
the isomorphism
$$H^i(S^{n-3} U^*
\otimes \lbrace \begin{smallmatrix}
1 & 2 &  3 & 2 & 1 \\
  &   &  1 &  &
\end{smallmatrix} \rbrace) =
H^i(S^{n-6}[ \begin{smallmatrix}
0 &  0 &  2 & 0 & 0 \\
  &   &   0 &  &
\end{smallmatrix} ]^*
\otimes \lbrace \begin{smallmatrix}
2 & 4 &  6 & 4 & 2 \\
  &   &  4 &  &
\end{smallmatrix} \rbrace).$$
That the latter cohomology vanishes for $i>0$
follows from Equation \ref{kos1:d_4(a_1)},
but with all the representations tensored with 
$\lbrace \begin{smallmatrix}
2 & 4 &  6 & 4 & 2 \\
  &   &  4 &  &
\end{smallmatrix} \rbrace$.
One needs to use Proposition \ref{sommers} (using the $A_5$ Levi subalgebra) 
on the initial term to show that 
$$H^i(S^{n-7}[ \begin{smallmatrix}
0 &  0 &  2 & 0 & 0 \\
  &   &   2 &  &
\end{smallmatrix} ]^*
\otimes \lbrace \begin{smallmatrix}
2 & 4 &  6 & 4 & 2 \\
  &   &  5 &  &
\end{smallmatrix} \rbrace) =
H^i(S^{n-10}[ \begin{smallmatrix}
0 &  0 &  2 & 0 & 0 \\
  &   &   2 &  &
\end{smallmatrix} ]^*
\otimes \lbrace \begin{smallmatrix}
3 & 6 &  9 & 6 & 3 \\
  &   &  5 &  &
\end{smallmatrix} \rbrace),$$
which is zero for $i>0$ by Proposition \ref{broer}.

Therefore, we have shown that there is an exact sequence 
$$0 \to 
H^0(S^{n-6}[ \begin{smallmatrix}
0 &  0 &  2 & 0 & 0 \\
  &   &   0 &  &
\end{smallmatrix} ]^*
\otimes \lbrace \begin{smallmatrix}
2 & 4 &  6 & 4 & 2 \\
  &   &  4 &  &
\end{smallmatrix} \rbrace)
\to 
H^0(S^{n} U^*) 
\to 
H^0(S^n[ \begin{smallmatrix}
0 &  1 &  0 & 1 & 0 \\
  &   &   1 &  &
\end{smallmatrix} ]^*) 
\to 
0.$$

\section{ $2A_2 + A_1$}

This orbit has normal closure.  We use the normality of 
the closure of $D_4(a_1)$ to prove the result.

Let $U_1 =  [ \begin{smallmatrix}
1  &  0 & 1 & 0 & 1 \\
  &   &   0 &  &
\end{smallmatrix} ]$.
Now $U_1$ is a sum of root spaces of $\g$ and 
we let $U_2$ be the subspace of $U_1$ of codimension two obtained  
by omitting the two root spaces for the roots
$\lbrace \begin{smallmatrix}
-1 & -1 &  -1 & 0 & 0 \\
  &   &  0 &  &
\end{smallmatrix} \rbrace$ and
$\lbrace \begin{smallmatrix}
0 & 0 &  -1 & -1 & -1 \\
  &   &  0 &  &
\end{smallmatrix} \rbrace$.  An application of Koszul and Proposition \ref{demazure}
shows that $S^n U_1^*$ 
and $S^n U_2^*$ have the same cohomology with respect to $G/B$.
Now let $U$ be the subspace of $\uni$ obtained by adding the root space of
the root 
$\lambda = \lbrace \begin{smallmatrix}
0 & -1 &  -1 & -1 & 0 \\
  &   &  -1 &  &
\end{smallmatrix} \rbrace$ to $U_2$.
Since $U$ is stable under $P_{\alpha_3}$ and 
$m = \langle -\lambda, \alpha_{3}\postcheck \rangle = -1$,
the cohomology of $S^n U^*$ and $S^n U_2^*$ coincide on $G/B$.
Consequently, $G \cdot U = G \cdot U_1$ and the latter 
equals the closure of $2A_2 + A_1$ as it arises from the weighted Dynkin diagram
for $2A_2 + A_1$.

We can prove normality by studying the short exact sequence
$$
0
\to
U
\to
[\begin{smallmatrix}
0 & 0 &  2 & 0 & 0 \\
  &   &  0 &  &
\end{smallmatrix}]
\to
V
\to
0
$$
and taking the Koszul resolution of its dual (there are eleven terms).

After considerable use of Proposition \ref{demazure}, 
it is possible to show that six of the nine initial terms 
of the resolution have total vanishing cohomology.  
The only possible non-zero contributions to cohomology occur for $n-9$ and $n-6$
and $n-3$
(the first, fourth, and seventh terms of the resolution).
All the other terms can be filtered so that the quotient line bundles have total
vanishing after one or more applications of Proposition \ref{demazure}.  
For example, in the sixth term of the resolution
a line bundle with weight 
$\lbrace \begin{smallmatrix}
0 & 2  & 4  & 2 & 0 \\
  &   & 2  &  &
\end{smallmatrix} \rbrace$ arises.
Using Proposition \ref{demazure}
we can replace this weight with the weight
$\lbrace \begin{smallmatrix}
0 & 2  & 4  & 2 & 1 \\
  &   & 2 &  &
\end{smallmatrix} \rbrace$ (if we shift cohomology degrees by one).
Then another application shows this weight will have total vanishing
cohomology since its pairing with $\al_4 \postcheck$ is $-1$.  

For the $n-3$ term, we need to proceed as in Lemma \ref{a1_thrice}.  The details are 
different, but the idea is the same. 
We study the bundle on the flag
variety of $A_2 + A_2 + A_1$ arising from the standard representation on each 
factor; this is an 18-dimensional representation which we denote by $U_3$.  
We consider the nine-dimensional $B$-subrepresentation $U'_3$
with weights corresponding to the weights which arise in $V^*$.
We can show that $\wedge^3 U'_3$ has Euler characteristic equal to zero 
and has cohomology (at worst) 
in degrees $2$ and $3$, where the cohomology is a sum of the two-dimensional
irreducible representation of $A_2 + A_2 + A_1$ (trivial on the first two factors, and
standard on the third).  But $S^3 U_3$ does not contain this representation and so
$\wedge^3 U'_3$ has no cohomology at all.  A spectral sequence argument 
as in \cite{demazure:bott-simple} yields the result 
we need in $E_6$.

For the $n-6$ term, we also need something akin to Lemma \ref{a1_thrice}. 
But first we consider the $B$-subrepresentation $Q$ of $\wedge^6 V^*$ containing
all $T$-weight spaces with weights $\lambda$ such that $\langle \lambda, 
\al_6 \postcheck \rangle$ is $-4$ or $-6$ (note that $Q$ extends to a representation
of the parabolic with Levi factor $A_2 + A_2 + A_1$).  We calculate that
$Q$ yields a bundle on $G/B$ with cohomology equal to the cohomology of the line bundle 
$\lbrace \begin{smallmatrix}
1 & 2  & 6  & 2 & 1 \\
  &   & 1  &  &
\end{smallmatrix} \rbrace$. 

We now want to 
show that the quotient of $\wedge^6 V^*$ by $Q$ has total vanishing cohomology.
The corresponding quotient of $\wedge^6 U'_3$ has Euler characteristic zero and
cohomology in (at worst) degrees $5$ and $6$, where the cohomology is a sum of
some number of copies of the trivial representation of $A_2 + A_2 + A_1$.  
But $S^6 U_3$ does not contain the trivial representation and so
there is no cohomology in this quotient of $\wedge^6 U'_3$.  
A spectral sequence argument yields that the quotient of
$\wedge^6 V^*$ by $Q$ has total vanishing.
We can thus conclude 
$$
H^{i+5}(S^{n-6} [ \begin{smallmatrix}
0 &  0 &  2 & 0 & 0 \\
  &   &   0 &  &
\end{smallmatrix} ]^*  \otimes 
\wedge^{6} V^*)
=
H^{i+5}(S^{n-6} [ \begin{smallmatrix}
0 &  0 &  2 & 0 & 0 \\
  &   &   0 &  &
\end{smallmatrix} ]^*  \otimes \lbrace \begin{smallmatrix}
1 & 2  & 6  & 2 & 1 \\
  &   & 1  &  &
\end{smallmatrix} \rbrace)
=
H^{i}(S^{n-6} [ \begin{smallmatrix}
0 &  0 &  2 & 0 & 0 \\
  &   &   0 &  &
\end{smallmatrix} ]^*  \otimes \lbrace \begin{smallmatrix}
2 & 4 &  6 & 4 & 2 \\
  &   &  4 &  &
\end{smallmatrix} \rbrace).$$

Finally,
$$H^{i+7}(S^{n-9} [ \begin{smallmatrix}
0 &  0 &  2 & 0 & 0 \\
  &   &   0 &  &
\end{smallmatrix} ]^*  \otimes 
\wedge^9 V^* )
=
H^{i+7}(S^{n-9} [ \begin{smallmatrix}
0 &  0 &  2 & 0 & 0 \\
  &   &   0 &  &
\end{smallmatrix} ]^*  \otimes \lbrace \begin{smallmatrix}
1 & 4 &  9 & 4 & 1 \\
  &   &  3 &  &
\end{smallmatrix} \rbrace)
=
H^i(S^{n-9} [ \begin{smallmatrix}
0 &  0 &  2 & 0 & 0 \\
  &   &   0 &  &
\end{smallmatrix} ]^*  \otimes \lbrace \begin{smallmatrix}
3 & 6 &  9 & 6 & 3 \\
  &   &  5 &  &
\end{smallmatrix} \rbrace).$$
Both of these line bundles can be shown to have cohomology which vanishes for $i>0$ by  
using Equation \ref{kos1:d_4(a_1)} with all the representations tensored by 
the appropriate weight
(we noted this for the first cohomology group in the previous section).
Thus we have the exact sequence
\begin{equation} \label{small:3}
0 \to
H^0(S^{n-9} [ \begin{smallmatrix}
0 &  0 &  2 & 0 & 0 \\
  &   &   0 &  &
\end{smallmatrix} ]^*  \otimes \lbrace \begin{smallmatrix}
3 & 6 &  9 & 6 & 3 \\
  &   &  5 &  &
\end{smallmatrix} \rbrace) \to
H^{0}(S^{n-6} [ \begin{smallmatrix}
0 &  0 &  2 & 0 & 0 \\
  &   &   0 &  &
\end{smallmatrix} ]^*  \otimes \lbrace \begin{smallmatrix}
2 & 4 &  6 & 4 & 2 \\
  &   &  4 &  &
\end{smallmatrix} \rbrace) \to
\end{equation} 
$$H^{0}(S^{n} [ \begin{smallmatrix}
0 &  0 &  2 & 0 & 0 \\
  &   &   0 &  &
\end{smallmatrix} ]^*) \to
H^0(S^n U^*) \to 0,$$
proving normality.


\section{ $A_2 + 2A_1$}

We prove the normality of the closure of $\orbit$ by
using the normality of the closure of $2A_2+A_1$.
The latter is the $G$-saturation of 
$U_1=[ \begin{smallmatrix}
1  &  0 & 1 & 0 & 1 \\
  &   &   0 &  &
\end{smallmatrix} ]$.
Consider the subspace $U$ of $U_1$ of codimension two obtained by
omitting the root spaces corresponding to the roots
$\lbrace \begin{smallmatrix}
-1 & -1 &  -1 & 0 & 0 \\
  &   &  0 &  &
\end{smallmatrix} \rbrace$ and
$\lbrace \begin{smallmatrix}
-1 & -1 &  -1 & 0 & 0 \\
  &   &  -1 &  &
\end{smallmatrix} \rbrace$.
It is possible to show that $G \cdot U$ equals
the closure of $\orbit$.  This is done by showing that there is a 
sequence of $B$-stable subspaces of $\uni$ beginning with 
$U$ and ending with $[ \begin{smallmatrix}
0  &  0 & 0 & 2 & 0 \\
  &   &   0 &  &
\end{smallmatrix} ]$ and such that each step in the sequence yields 
the needed isomorphism of cohomology of symmetric powers of dual spaces.
The subspace $[ \begin{smallmatrix}
0  &  0 & 0 & 2 & 0 \\
  &   &   0 &  &
\end{smallmatrix} ]$
has saturation equal to the closure of $\orbit$.



Thus the relevant short exact sequence is 
$$
0
\to
U
\to
[ \begin{smallmatrix}
1  &  0 & 1 & 0 & 1 \\
  &   &   0 &  &
\end{smallmatrix} ]
\to
V
\to
0
,$$ with Koszul resolution of its dual equal to 

\begin{equation}\label{a2+2a1}
0 \to
S^{n-2} [ \begin{smallmatrix}
1  &  0 & 1 & 0 & 1 \\
  &   &   0 &  &
\end{smallmatrix} ]^*  \otimes
\wedge^2 V^*
\to
S^{n-1} [ \begin{smallmatrix}
1  &  0 & 1 & 0 & 1 \\
  &   &   0 &  &
\end{smallmatrix} ]^*  \otimes
V^*
\to
S^{n} [ \begin{smallmatrix}
1  &  0 & 1 & 0 & 1 \\
  &   &   0 &  &
\end{smallmatrix} ]^*
\to
S^{n} U^*
\to
0
.
\end{equation}
It is easy to see that 
$H^i(S^{n-1} [ \begin{smallmatrix}
1  &  0 & 1 & 0 & 1 \\
  &   &   0 &  &
\end{smallmatrix} ]^*  \otimes 
V^*)
=0$ for all $i,n$ by two applications of Proposition \ref{demazure},
and that 
$$
H^{i+1}(S^{n-2} [ \begin{smallmatrix}
1  &  0 & 1 & 0 & 1 \\
  &   &   0 &  &
\end{smallmatrix} ]^*  \otimes
\wedge^2 V^*)
=
H^{i+1}(S^{n-2} [ \begin{smallmatrix}
1  &  0 & 1 & 0 & 1 \\
  &   &   0 &  &
\end{smallmatrix} ]^*  \otimes \lbrace \begin{smallmatrix}
2 &  2 &  2 & 0  & 0 \\
  &   &  1 &  &
\end{smallmatrix} \rbrace)=
H^{i}(S^{n-2} [ \begin{smallmatrix}
1  &  0 & 1 & 0 & 1 \\
  &   &   0 &  &
\end{smallmatrix} ]^*  \otimes \lbrace \begin{smallmatrix}
2 &  2 &  2 & 1 & 0 \\
  &   &  1 &  &
\end{smallmatrix} \rbrace).$$
Using Equation \ref{a2+2a1} again, but tensoring each term with
$\lbrace \begin{smallmatrix}
0 &  1 &  2 & 2  & 2 \\
  &   &  1 &  &
\end{smallmatrix} \rbrace$,
we can show that the latter cohomology group coincides with
$$H^i(S^{n-4} [ \begin{smallmatrix}
1  &  0 & 1 & 0 & 1 \\
  &   &   0 &  &
\end{smallmatrix} ]^*  \otimes \lbrace \begin{smallmatrix}
2 &  3 &  4 & 3  & 2 \\
  &   &  2 &  &
\end{smallmatrix} \rbrace)$$ as long as we can
first show that 
$$H^i(S^{n-2} U^* \otimes \lbrace \begin{smallmatrix}
0 &  1 &  2 & 2  & 2 \\
  &   &  1 &  &
\end{smallmatrix} \rbrace)$$ vanishes for all $i \geq 0$.
This amounts to using the sequence of subspaces which
connects $U$ to $[ \begin{smallmatrix}
0  &  0 & 0 & 2 & 0 \\
  &   &   0 &  &
\end{smallmatrix} ]$ (see the first paragraph above) 
and transferring the problem to one using 
this latter subspace. 

Finally, $$H^i(S^{n-4} [ \begin{smallmatrix}
1  &  0 & 1 & 0 & 1 \\
  &   &   0 &  &
\end{smallmatrix} ]^*  \otimes \lbrace \begin{smallmatrix}
2 &  3 &  4 & 3  & 2 \\
  &   &  2 &  &
\end{smallmatrix} \rbrace)$$ has higher vanishing by Proposition \ref{broer}, our first
application of this proposition in its full generality
(here, $\omega = 
\lbrace \begin{smallmatrix}
-2 &  -5 &  -8 & -5  & -2 \\
  &   &  -4 &  &
\end{smallmatrix} \rbrace$
and
$\lambda =
\lbrace \begin{smallmatrix}
4 &  8 &  12 & 8  & 4 \\
  &   &  6 &  &
\end{smallmatrix} \rbrace$).  Hence the exact sequence
$$0 \to 
H^0(S^{n-4} [ \begin{smallmatrix}
1  &  0 & 1 & 0 & 1 \\
  &   &   0 &  &
\end{smallmatrix} ]^*  \otimes \lbrace \begin{smallmatrix}
2 &  3 &  4 & 3  & 2 \\
  &   &  2 &  &
\end{smallmatrix} \rbrace) \to
H^{0}(S^{n} [ \begin{smallmatrix}
1  &  0 & 1 & 0 & 1 \\
  &   &   0 &  &
\end{smallmatrix} ]^* ) \to
H^{0}(S^n U^*) \to 0,$$
proving normality.





\section{ $A_2$}

We deduce the normality of $\bar{\orbit}$ from the normality 
of the closure of $A_2 + 2A_1$.
First, we claim that $\bar{\orbit}$ is the $G$-saturation 
of the subspace 
$$U_2 = [ \begin{smallmatrix}
0  &  0 & 0 & 0 & 0 \\
  &   &   2 &  &
\end{smallmatrix}] 
\cap
[ \begin{smallmatrix}
0  &  1 & 0 & 1 & 0 \\
 &   &   0 &  &
\end{smallmatrix} ].$$
This is proved by using the Koszul resolution of the dual of 
\begin{equation}\label{a2}
0
\to
U_2 
\to
[ \begin{smallmatrix}
0  &  0 & 0 & 0 & 0 \\
 &   &   2 &  &
\end{smallmatrix} ] 
\to
V_1
\to
0
\end{equation} 
to show that
$$H^i(S^n U_2^*) = H^i(S^n[ \begin{smallmatrix}
0  &  0 & 0 & 0 & 0 \\
  &   &   2 &  &
\end{smallmatrix}]^*)$$
for all $i,n$.  We omit the details.

Next, we study the Koszul resolution of the dual of 
$$0
\to
U_2
\to
[ \begin{smallmatrix}
0  &  1 & 0 & 1 & 0 \\
 &   &   0 &  &
\end{smallmatrix} ] 
\to
V_2
\to
0.$$
The weights of $V^*_2$ are
$\lbrace \begin{smallmatrix}
0 &  1 &  1 & 1  & 0 \\
  &   &  0 &  &
\end{smallmatrix} \rbrace$,
$\lbrace \begin{smallmatrix}
1 &  1 &  1 & 1  & 0 \\
  &   &  0 &  &
\end{smallmatrix} \rbrace$,
$\lbrace \begin{smallmatrix}
0 &  1 &  1 & 1  & 1 \\
  &   &  0 &  &
\end{smallmatrix} \rbrace$, and 
$\lbrace \begin{smallmatrix}
1 &  1 &  1 & 1  & 1 \\
  &   &  0 &  &
\end{smallmatrix} \rbrace$,
and we find that the cohomology of the kernel $C$ of the map
$$S^n [\begin{smallmatrix}
0  &  1 & 0 & 1 & 0 \\
  &   &   0 &  &
\end{smallmatrix} ]^*
\to
S^n U_2^*$$ 
satisfies the long exact sequence
$$\dots \to H^i(S^{n-4} [ \begin{smallmatrix}
0  &  1 & 0 & 1 & 0 \\
  &   &   0 &  &
\end{smallmatrix} ]^*  \otimes \lbrace \begin{smallmatrix}
2 &  4 &  6 & 4  & 2 \\
  &   &  3 &  &
\end{smallmatrix} \rbrace) 
\to
H^i(S^{n-3} [ \begin{smallmatrix}
0  &  1 & 0 & 1 & 0 \\
  &   &   0 &  &
\end{smallmatrix} ]^*  \otimes \lbrace \begin{smallmatrix}
2 &  3 &  4 & 3  & 2 \\
  &   &  2 &  &
\end{smallmatrix} \rbrace)
\to 
H^i(C) 
\to 
\dots$$ 

If we can show that the first term above vanishes for $i \geq 2$ and
the second term vanishes for $i \geq 1$, this will be sufficient to 
deduce that $H^i(C)=0$ for $i \geq 1$, and normality will follow.  
We sketch our argument.

Let $\lambda_1 = \lbrace \begin{smallmatrix}
2 &  4 &  6 & 4  & 2 \\
  &   &  3 &  &
\end{smallmatrix} \rbrace$
and $\lambda_2 = \lbrace \begin{smallmatrix}
2 &  3 &  4 & 3  & 2 \\
  &   &  2 &  &
\end{smallmatrix} \rbrace$.
Let $U_3 = [ \begin{smallmatrix}
0  &  1 & 0 & 1 & 0 \\
  &   &   0 &  &
\end{smallmatrix} ] \cap
[ \begin{smallmatrix}
1  &  0 & 1 & 0 & 1 \\
  &   &   0 &  &
\end{smallmatrix} ]$. 

The Koszul resolution coming from the inclusion of
$U_3$ into $[ \begin{smallmatrix}
0  &  1 & 0 & 1 & 0 \\
  &   &   0 &  &
\end{smallmatrix}]$ 
is 
$$0
\to
S^{n-2} [ \begin{smallmatrix}
0  &  1 & 0 & 1 & 0 \\
  &   &   0 &  &
\end{smallmatrix} ]^*  \otimes \lbrace \begin{smallmatrix}
0 &  2 &  2 & 2  & 0 \\
  &   &  1 &  &
\end{smallmatrix} \rbrace 
\to
S^{n-1} [ \begin{smallmatrix}
0  &  1 & 0 & 1 & 0 \\
  &   &   0 &  &
\end{smallmatrix} ]^*  \otimes V^*_3
\to
S^n [ \begin{smallmatrix}
0  &  1 & 0 & 1 & 0 \\
  &   &   0 &  &
\end{smallmatrix} ]^*
\to
S^n U^*_3
\to 0.$$
 
Tensoring this equation with $\lambda_2$ and then taking cohomology,
we find that the first two terms have total vanishing cohomology
since the weights of $V^*_3$ are 
$\lbrace \begin{smallmatrix}
0 &  1 &  1 & 1  & 0 \\
  &   &  1 &  &
\end{smallmatrix} \rbrace$ and
$\lbrace \begin{smallmatrix}
0 &  1 &  1 & 1  & 0 \\
  &   &  0 &  &
\end{smallmatrix} \rbrace$. 
On the other hand, tensoring with $\lambda_1$, the second term has vanishing
cohomology, but the cohomology of the first term in degree $i$ coincides 
with the cohomology
$H^{i-2}(S^{n-2} [ \begin{smallmatrix}
0  &  1 & 0 & 1 & 0 \\
  &   &   0 &  &
\end{smallmatrix} ]^*  \otimes \lbrace \begin{smallmatrix}
3 &  6 &  8 & 6  & 3 \\
  &   &  4 &  &
\end{smallmatrix} \rbrace)$.
Since the latter does vanish for $i>2$ by Proposition \ref{broer}, we have
reduced our question to showing that
$H^i(S^{n-4} U_3^*  \otimes \lambda_1)=0$ for $i \geq 2$
and 
$H^i(S^{n-3} U_3^* \otimes \lambda_2)=0$ for $i \geq 1$.

We can prove these two results by studying 
the Koszul resolution of the dual of 
\begin{equation} \label{upit}
0
\to
U_3
\to
[ \begin{smallmatrix}
1  &  0 & 1 & 0 & 1 \\
  &   &   0 &  &
\end{smallmatrix} ]
\to
V_4
\to
0
.
\end{equation}

The weights of $V^*_4$ are 
$\lbrace \begin{smallmatrix}
1 &  1 &  1 & 0  & 0 \\
  &   &  0 &  &
\end{smallmatrix} \rbrace$,
$\lbrace \begin{smallmatrix}
1 &  1 &  1 & 0  & 0 \\
  &   &  1 &  &
\end{smallmatrix} \rbrace$,
$\lbrace \begin{smallmatrix}
0 &  0 &  1 & 1  & 1 \\
  &   &  0 &  &
\end{smallmatrix} \rbrace$, and
$\lbrace \begin{smallmatrix}
0 &  0 &  1 & 1  & 1 \\
  &   &  1 &  &
\end{smallmatrix} \rbrace$.

Tensoring the Koszul resolution of the dual of (\ref{upit}) with $\lambda_2$, we find that we must study 
\begin{equation} \label{key_coho}
H^i(S^n [ \begin{smallmatrix}
1  &  0 & 1 & 0 & 1 \\
  &   &   0 &  &
\end{smallmatrix} ]^* \otimes \mu)
\end{equation} 
where $\mu$ is one of the four weights $\lambda_2$, $\lbrace \begin{smallmatrix}
2 &  4 & 6  & 5  & 4 \\
  &   &  3 &  &
\end{smallmatrix} \rbrace$,
$\lbrace \begin{smallmatrix}
4 &  5 &  6 & 4  & 2 \\
  &   &  3 &  &
\end{smallmatrix} \rbrace$, and
$2 \lambda_2$.
Although we do not need it in its full strength, 
we find that all four have vanishing for $i>0$.
For $\mu= \lambda_2$ or $2 \lambda_2$, the cohomology vanishes for $i>0$ 
by Proposition \ref{broer}.
For the other two weights, we must use Equation \ref{a2+2a1} from the previous section.
Tensoring that equation with 
$\mu = \lbrace \begin{smallmatrix}
2 &  4 & 6  & 5  & 4 \\
  &   &  3 &  &
\end{smallmatrix} \rbrace$, 
we must study the term
$H^i(S^n [ \begin{smallmatrix}
1  &  0 & 1 & 0 & 1 \\
  &   &   0 &  &
\end{smallmatrix} ]^* \otimes 2 \lambda_2)$, which we just said vanished for $i>0$,
and the term 
$H^i(S^n U^* \otimes \mu)$.
Using the method from the previous section, we can show
that the latter coincides with
$H^i(S^n [ \begin{smallmatrix}
0  &  0 & 0 & 2 & 0 \\
  &   &   0 &  &
\end{smallmatrix} ]^* \otimes \mu)$.
This vanishes for $i>0$:
we study the inclusion of 
$[ \begin{smallmatrix}
0  &  0 & 0 & 2 & 0 \\
  &   &   0 &  &
\end{smallmatrix} ]$ into $[ \begin{smallmatrix}
0  &  0 & 0 & 2 & 2 \\
  &   &   0 &  &
\end{smallmatrix} ]$.
We must show that 
$H^i(S^n [ \begin{smallmatrix}
0  &  0 & 0 & 2 & 2 \\
  &   &   0 &  &
\end{smallmatrix}]^* \otimes 
\lbrace \begin{smallmatrix}
2 &  4 &  6 & 5  & 5 \\
  &   &  3 &  &
\end{smallmatrix} \rbrace)$ vanishes for $i>0$.
This follows since it equals 
$H^i(S^n [ \begin{smallmatrix}
0  &  0 & 0 & 0 & 2 \\
  &   &   2 &  &
\end{smallmatrix}]^* \otimes 
\lbrace \begin{smallmatrix}
3 &  6 &  9 & 7  & 5 \\
  &   &  5 &  &
\end{smallmatrix} \rbrace)$
by a result analogous to Proposition \ref{sommers} and we can now
invoke Proposition \ref{broer}.  The final $\mu$ is handled in a symmetric fashion to this one.

The situation for $\lambda_1$ 
requires us to study the cohomologies in 
Equation \ref{key_coho} for $\mu$ equal to $\nu + \lbrace \begin{smallmatrix}
0 &  1 &  2 & 1  & 0 \\
  &   &  1 &  &
\end{smallmatrix} \rbrace$ where $\nu$ is one of the four 
weights listed for $\lambda_2$.

We can analyze these weights as follows.
Let $V_5$ equal $[ \begin{smallmatrix}
1  &  0 & 1 & 0 & 1 \\
  &   &   0 &  &
\end{smallmatrix} ]$ but omitting the root space 
$\lbrace \begin{smallmatrix}
0 &  -1 &  -2 & -1  & 0 \\
  &   &  -1 &  &
\end{smallmatrix} \rbrace$.
Then the cohomology results we need will follow if
we can show that $H^i(S^n V^*_5 \otimes \nu) =0 $ for $i>0$.
This follows by studying the inclusion of $V_5$ into
$[ \begin{smallmatrix}
2  &  0 & 0 & 0 & 2 \\
  &   &   0 &  &
\end{smallmatrix} ]$.  We omit these details.
  
\section{ $3A_1$}

The closure of the orbit $\orbit$ of type $3A_1$ 
is the $G$-saturation of 
$[\begin{smallmatrix}
0 & 0 &  1 & 0 & 0 \\
  &   &  0 &  &
\end{smallmatrix}]$.
Let $U_4$ be the $B$-stable subspace of $\uni$ obtained from 
$[\begin{smallmatrix}
0 & 0 &  1 & 0 & 0 \\
  &   &  0 &  &
\end{smallmatrix}]$
by omitting the root space for 
$\lbrace \begin{smallmatrix}
0 &  -1 &  -2 & -1  & 0 \\
  &   &  -1 &  &
\end{smallmatrix} \rbrace$ and 
adding the root space for 
$\lbrace \begin{smallmatrix}
-1 &  -1 &  -1 & -1  &  -1\\
  &   &  -1 &  &
\end{smallmatrix} \rbrace$.
Then the $G$-saturation of $U_4$ equals $\bar{\orbit}$; we omit the details.  

Next let $U_2 = [ \begin{smallmatrix}
0  &  0 & 0 & 0 & 0 \\
  &   &   2 &  &
\end{smallmatrix}] 
\cap
[ \begin{smallmatrix}
0  &  1 & 0 & 1 & 0 \\
 &   &   0 &  &
\end{smallmatrix} ].$  As we noted in the previous section,
there is an isomorphism  
$$H^i(S^n U_2^*) = H^i(S^n[ \begin{smallmatrix}
0  &  0 & 0 & 0 & 0 \\
  &   &   2 &  &
\end{smallmatrix}]^*)$$
induced by the inclusion of $U_2$ into $[ \begin{smallmatrix}
0  &  0 & 0 & 0 & 0 \\
  &   &   2 &  &
\end{smallmatrix}]$.

Consider the short exact sequence 
$$
0
\to
U_4
\to
U_2
\to
V
\to
0.
$$
The analysis of the Koszul resolution of its dual leads to the exact
sequence
$$
\dots \to H^i(S^{n-4} U_2^*
\otimes \lbrace \begin{smallmatrix}
2 &  4 &  6 & 4  & 2 \\
  &   &  4 &  &
\end{smallmatrix} \rbrace) 
\to 
H^i(S^{n} U_2^*)
\to 
H^i(S^{n} U_4^*)
\to \dots
$$
The proof uses the fact that $U_2$ is stable for $P_{\al_1}, P_{\al_3}, P_{\al_5}$
and the representation $V^*$ restricts to the representation in Lemma \ref{a1_thrice}.
The result follows as in Section \ref{section_a5}.

Equation \ref{a2} tensored with
$\lbrace \begin{smallmatrix}
2 &  4 &  6 & 4  & 2 \\
  &   &  4 &  &
\end{smallmatrix} \rbrace$ leads to the isomorphism
$$H^i(S^{n-4} U^*_2
\otimes \lbrace \begin{smallmatrix}
2 &  4 &  6 & 4  & 2 \\
  &   &  4 &  &
\end{smallmatrix} \rbrace) =
H^i(S^{n-4}
[\begin{smallmatrix}
0 & 0 &  0 & 0 & 0 \\
  &   &  2 &  &
\end{smallmatrix}]^*
\otimes \lbrace \begin{smallmatrix}
2 &  4 &  6 & 4  & 2 \\
  &   &  4 &  &
\end{smallmatrix} \rbrace).$$ 
As the latter vanishes for $i>0$, we have the desired surjectivity
(induced by inclusions) of
$H^0(S^n U^*_2)$ $= H^0(S^n[ \begin{smallmatrix}
0  &  0 & 0 & 0 & 0 \\
  &   &   2 &  &
\end{smallmatrix}]^*)$ onto 
$H^0(S^n U_4^*)$ 
with kernel isomorphic to 
$H^0(S^{n-4}
[\begin{smallmatrix}
0 & 0 &  0 & 0 & 0 \\
  &   &  2 &  &
\end{smallmatrix}]^*
\otimes \lbrace \begin{smallmatrix}
2 &  4 &  6 & 4  & 2 \\
  &   &  4 &  &
\end{smallmatrix} \rbrace)$,
proving normality.

\section{ $2A_1$}

This orbit is already known to have normal closure by Hesselink \cite{hesselink:normal}.
In any event we can also prove it by showing that there is an exact sequence
$$
0
\to
H^0(S^{n-3}[\begin{smallmatrix}
0 & 0 &  1 & 0 & 0 \\
  &   &  0 &  &
\end{smallmatrix}]^* \otimes \lbrace \begin{smallmatrix}
2 &  4 &  6 & 4  & 2 \\
  &   &  3 &  &
\end{smallmatrix} \rbrace)
\to
H^0(S^n[\begin{smallmatrix}
0 & 0 &  1 & 0 & 0 \\
  &   &  0 &  &
\end{smallmatrix}]^*)
\to
H^0(S^n[\begin{smallmatrix}
1 & 0 &  0 & 0 & 1 \\
  &   &  0 &  &
\end{smallmatrix}]^*)
\to
0.$$

\section{ $A_1$}

This has normal closure by Hesselink \cite{hesselink:normal} or 
Vinberg-Popov \cite{vinberg-popov}.

\section{The non-normal nilpotent varieties}

The orbits $A_4$, $A_3+A_1$, $A_3$, $2A_2$, $A_2+A_1$
all have non-normal closure.
The easiest way to see this is to show that the induced map 
$H^0(S^n \uni^*) \to H^0(S^n V^*)$ is not surjective, where
$V$ is as in Case 2 of Proposition \ref{broer} with $G$-saturation
the desired orbit closure.

One calculates by hand (invoking McGovern \cite{mcgovern:regfn}) that
the adjoint representation has non-zero multiplicity in 
$H^0(S^n V^*)$ for $n=3$ for $A_4$,$A_3+A_1$, and $A_3$; and 
for $n=2$ for $2A_2$ and $A_2+A_1$.
On the other hand, by Kostant \cite{kostant:poly}, the adjoint
representation has non-zero multiplicity in $H^0(S^n \uni^*)$ 
only when $n=1, 4, 5, 7, 8, 11$ (the exponents of $E_6$).

See \cite{kraft:normalG2} and \cite{broer:normalF4} for a survey of 
techniques to show that an orbit closure is not normal.

\section{Proof of Theorem \ref{small_conj}}

For each of the pairs in the theorem we showed above that the functions
on the first orbit of degree $n$ 
are a quotient of the functions on the second orbit 
(or a cover of it) of degree $n$ and computed the kernel.  
We have
\begin{enumerate}
\item For $(A_5, E_6(a_3))$, the kernel is  
$H^0(S^{n-10} [ \begin{smallmatrix}
0 &  2 &  0 & 2 & 0 \\
  &   &   2 &  &
\end{smallmatrix} ]^*  \otimes \lbrace \begin{smallmatrix}
2 & 4 &  6 & 4 & 2 \\
  &   &  4 &  & 
\end{smallmatrix} \rbrace)$.
\item For $(2A_2+A_1, D_4(a_1))$,
the kernel is a quotient of 
$H^0(S^{n-6}[ \begin{smallmatrix}
0 &  0 &  2 & 0 & 0 \\
  &   &   0 &  &
\end{smallmatrix} ]^*
\otimes \lbrace \begin{smallmatrix}
2 & 4 &  6 & 4 & 2 \\
  &   &  4 &  &
\end{smallmatrix} \rbrace)$.
\item For $(A_3+A_1, \tilde{D_4}(a_1))$, 
the kernel is $H^0(S^{n-6}[ \begin{smallmatrix}
0 &  0 &  2 & 0 & 0 \\
  &   &   0 &  &
\end{smallmatrix} ]^*
\otimes \lbrace \begin{smallmatrix}
2 & 4 &  6 & 4 & 2 \\
  &   &  4 &  &
\end{smallmatrix} \rbrace)$.
\item For $(3A_1, A_2)$, the kernel is 
$H^0(S^{n-4}
 [\begin{smallmatrix}
 0 & 0 &  0 & 0 & 0 \\
   &   &  2 &  &
 \end{smallmatrix}]^*
 \otimes \lbrace \begin{smallmatrix}
 2 &  4 &  6 & 4  & 2 \\
   &   &  4 &  &
\end{smallmatrix} \rbrace)$.
\end{enumerate}

Since the higher cohomologies of these bundles vanish, we can
compute the multiplicity of any finite-dimensional representation in $H^0(-)$ by using
the Bott-Borel-Weil Theorem.
The fact that 
$\lbrace \begin{smallmatrix}
 2 &  4 &  6 & 4  & 2 \\
   &   &  4 &  &
\end{smallmatrix} \rbrace$ is twice a root (the highest root
of $E_6$) implies immediately that no small representation
has non-zero multiplicity in $H^0(-)$, proving the theorem. 

\section{Conclusions}

We can use the same techniques to prove that many orbit closures in $E_7$
and $E_8$ are normal.   However, since we were not able to resolve the 
picture completely in those types, 
we did not include those calculations in this paper.
We can also use these techniques to resolve the analog
of Theorem \ref{small_conj} in
types $G_2$ and partially in $E_7$ and $E_8$ in the same manner as we did here.

To extend these results to good positive characteristic one would
have to find a substitute for the use of Proposition \ref{broer}.  
Propositions \ref{demazure} and \ref{sommers}, however,
can be shown (in the same vein as \cite{thomsen:normal}) to carry over in the generality
that we used them here.

\bibliography{sommers1}
\bibliographystyle{pnaplain}
\end{document}